\documentclass[preprint,12pt,3p]{elsarticle}

\usepackage{amssymb}
\usepackage{amsmath,amsthm}
\usepackage{mathrsfs}
\usepackage{hyperref}
\usepackage{cleveref}
\usepackage[all,cmtip]{xy}
\usepackage{epsf, graphicx}
\usepackage{latexsym,amsfonts,amsbsy,amssymb}
\usepackage{geometry}
\usepackage{titletoc}
\usepackage{algorithm}
\usepackage{algorithmic}
\usepackage{amscd}

\usepackage{amssymb}

\newtheorem{theorem}{Theorem}
\newtheorem{lemma}[theorem]{Lemma}
\newtheorem{corollary}[theorem]{Corollary}
\newtheorem{notation}[theorem]{Notation}

\makeatletter
\def\ps@pprintTitle{%
\let\@oddhead\@empty
\let\@evenhead\@empty
\def\@oddfoot{\reset@font\hfil\thepage\hfil}
\let\@evenfoot\@oddfoot
}
\makeatother

\def\O{{\cal O}}

\begin{document}

\begin{frontmatter}

\title{Descent of equivalences for blocks with Klein four defect groups\tnoteref{label0}}
\tnotetext[label0]{The author gratefully acknowledges financial support by CSC (No. 202006770016) and a program of CCNU (No. 2020CXZZ074).}

\author[label1,label2]{Xin Huang}
\address[label1]{School of Mathematics and Statistics, Central China Normal University,
Wuhan, 430079, P.R. China}
\address[label2]{Department of Mathematics, City, University of London, London, EC1V 0HB, United Kingdom}


\ead{xinhuang@mails.ccnu.edu.cn}




\begin{abstract}
We show that the splendid Rickard complexes for blocks with Klein four defect groups constructed by Rickard and Linckelmann descend to non-split fields. As a corollary, Navarro's refinement of the Alperin-McKay conjecture holds for blocks with a Klein four defect group. We also prove that splendid Morita equivalences between blocks and their Brauer correspondents (if exist) descend to non-split situations.
\end{abstract}


\end{frontmatter}




As noted in \cite{Kessar_Linckelmann}, categorical equivalences between block algebras of finite groups, such as Morita equivalences and Rickard equivalences,  induce character bijections which commute with the Galois groups of field extensions (see \cite[Theorem 1.3, 1.6 and Corollary 1.9]{Kessar_Linckelmann}). This is the motivation for realising known Morita and Rickard equivalences over non splitting fields. Kessar and Linckelmann showed that Rouquier's splendid Rickard complexes for blocks with cyclic defect groups descend to non-split fields (see \cite[Theorem 1.10]{Kessar_Linckelmann}). By using it, they recovered the result \cite[Theorem 3.4]{Navarro} of Navarro (see \cite[Corollary 1.11]{Kessar_Linckelmann}). Namely, the Alperin-McKay-Navarro conjecture (\cite[Conjecture B]{Navarro}) holds for blocks with cyclic defect groups.
Combining \cite[Theorem 1.10]{Kessar_Linckelmann} with the general descent results (\cite[Theorem 6.5]{Kessar_Linckelmann}), Kessar and Linckelmann proved that Brou\'{e}'s abelian defect group conjecture (originally stated over complete discrete valuation rings with splitting residue fields) holds for blocks with cyclic defect groups over arbitrary complete discrete valuation rings (see \cite[Theorem 1.12]{Kessar_Linckelmann}).

Throughout this note, $p$ is a prime number, $k\subseteq k'$ are fields of characteristic $p$; $\O\subseteq\O'$ are either complete discrete valuation rings of characteristic $0$ with $J(\O)\subseteq J(\O')$ and with residue fields $k,~k'$ respectively, or $\O=k,~\O'=k'$.

For a finite group $G$, a {\it block} of the group algebra $\O G$ is a primitive idempotent of the center of $\O G$. For a subgroup $H$ of $G$, let $(\O G)^H$ denote the set of $H$-fixed elements of $\O G$ under the conjugation action. If $H$ is a $p$-subgroup, the {\it Brauer map} is the $\O$-algebra homomorphism
${\rm Br}_H: (\O G)^H\to kC_G(H)$, $\sum_{g\in G}\alpha_gg\mapsto \sum_{g\in C_G(H)}\bar{\alpha}_gg,$
where $\bar{\alpha}_g$ denotes the image of $\alpha_g$ in $k$.
For a block $b$ of $G$, a {\it defect group} of $b$ is a maximal $p$-subgroup $P$ of $G$ such that ${\rm Br}_P(b)\neq 0$. By Brauer's first main theorem, there is a unique block $c$ of $\O N_G(P)$ with defect group $P$ such that ${\rm Br}_P(b) = {\rm Br}_P(c)$ and the map $b\mapsto c$ is a bijection between the set of blocks of $\O G$ with defect group $P$ and the set of blocks of $\O N_G(P)$ with defect group $P$. This bijection is known as the {\it Brauer correspondence}.



Inspired by \cite[Theorem 1.10]{Kessar_Linckelmann} and its proof, we prove the following result. For the definition of the splendid Rickard equivalences, we refer to \cite[Definition 9.7.5]{Linckelmann}. The proof uses \cite[Theorem 1.1]{CEKL}, which relies on the classification of finite simple groups.

\begin{theorem}\label{3.1}
Let $G$ be a finite group, $b$ a block of $\O'G$ having a Klein four group defect group $P$. Let $c$ be the block of $\O'N_G(P)$ corresponding to $b$ via the Brauer correspondence. Suppose that $b\in \O G$ and that $k'$ is a splitting field for all subgroups of $G$. Then the block algebras $\O Gb$ and $\O N_G(P)c$ are splendidly Rickard equivalent. More precisely, there is a splendid Rickard complex $X$ of $(\O N_G(P)c, \O Gb)$-bimodules such that $\O'\otimes_{\O} X$ is isomorphic to the complex in {\rm \cite[Theorem 12.4.1]{Linckelmann}}.
\end{theorem}

By \cite[Corollary 1.9]{Kessar_Linckelmann} and the remarks following it, Theorem \ref{3.1} has the following corollary.

\begin{corollary}
{\rm \cite[Conjecture B]{Navarro}} holds for blocks with a Klein four defect group.
\end{corollary}

Using \cite[Theorem 6.5]{Kessar_Linckelmann} and Theorem \ref{3.1}, we can prove the following result. As \cite[Theorem 1.12]{Kessar_Linckelmann} and \cite[Theorem A, B]{Boltje}, it may be viewed as evidence for the Brou\'{e}'s abelian defect group conjecture over arbitrary complete discrete valuation rings.

\begin{theorem}\label{main}
Let $G$ be a finite group, $b$ a block of $\O G$ having a Klein four defect group. Then $\O Gb$ is splendidly Rickard equivalent to its Brauer correspondent block algebra.
\end{theorem}

By using the same tool (i.e., Lemma \ref{key} below), we can also prove the following theorem.

\begin{theorem}\label{3}
Assume that $\O'$ is finitely generated as an $\O$-module. Let $G$ be a finite group, $b$ a block of $\O'G$ having a defect group $P$. Let $c$ be the block of $\O'N_G(P)$ corresponding to $b$ via the Brauer correspondence. Suppose that $b\in \O G$ and that $k'$ is a splitting field for all subgroups of $G$. Assume that the block algebras $\O'Gb$ and $\O'N_G(P)c$ are splendidly Morita equivalent. Then the block algebras $\O Gb$ and $\O N_G(P)c$ are splendidly Morita equivalent. More precisely, the bimodule $T'$ (resp. $T$) in Lemma \ref{key} induces a Morita equivalence between $\O'Gb$ and $\O'N_G(P)c$ (resp. between $\O Gb$ and $\O N_G(P)c$).
\end{theorem}

\begin{corollary}\label{corollary2}
Let $G$ be a finite group, $b$ a block of $\O'G$. Suppose that $k'$ is a finite splitting field for all subgroups of $G$ and that $\O'$ has characteristic 0 and contains a primitive $|G|$-th root of unity. Assume that the block algebras $\O'Gb$ is splendidly Morita equivalent to its Brauer correspondent block algebra.   Then {\rm \cite[Conjecture B]{Navarro}} holds for $b$.
\end{corollary}



Let $G$ be a finite group, $b$ a block of $\O G$ and $P$ a defect group of $b$. A primitive idempotent $i\in (\O Gb)^P$ satisfying ${\rm Br}_P(i)\neq 0$ is called a {\it source idempotent} of $b$. Then the interior $P$-algebra $i\O Gi$, with the structure homomorphism $P\to (i\O Gi)^\times$, $u\mapsto ui=iu$, is called a
{\it source algebra} of $b$. By \cite[Corollary 3.5]{Pointed}, the $(\O Gb,i\O Gi)$-bimodule $\O Gi$ and the $(i\O Gi,\O Gb)$-bimodule $i\O G$ induce a Morita equivalence between $\O Gb$ and $i\O Gi$.
Assume that $P$ is a Klein four group and $k$ is a splitting field for all subgroups of $G$. \cite[Theorem 1.1]{CEKL} tells us that, by using the classification of finite simple groups, every source algebra of $b$ is isomorphic, as an interior $P$-algebra, to either $\O P$, $\O A_4$ or $\O A_5b_0$, where $b_0$ is the principle block idempotent of $\O A_5$. As noted in \cite{CEKL}, here $\O A_4$ and the $\O A_5b_0$ are viewed as interior $P$-algebras via some identification of $P$ with the Sylow $2$-subgroups of $A_4$ and $A_5$. The interior $P$-algebra
structure is independent of the choice of an identification because any automorphism of a
Sylow $2$-subgroup of $A_4$ or $A_5$ extends to an automorphism of $A_4$ or $A_5$, respectively.

\begin{notation}\label{notation}
{\rm Let $G$ be a finite group, $b$ a block of $\O'G$ having a defect group $P$. Let $e$ be a block of $k'C_G(P)$ such that ${\rm Br}_P(b)e=e$. Let $\hat{e}$ be the block of $\O' C_G(P)$ that lifts the block $e$ of $k'C_G(P)$. Let $c$ be the block of $\O'N_G(P)$ corresponding to $b$ via the Brauer correspondence. Let $j\in\O' C_G(P)\hat{e}$ be a primitive idempotent. Then $j$ is a source idempotent of $\hat{e}$ as a block of $\O' N_G(P,e)$, and also a source idempotent of the block $c$ of $\O' N_G(P)$ (see \cite[Theorem 6.2.6 (iv)]{Linckelmann}). Moreover, $j\O' N_G(P,e)j=j\O'N_G(P)j$ (see \cite[Theorem 6.8.3]{Linckelmann}). By \cite[Proposition 4.10]{FP} or \cite[Proposition 6.7.4]{Linckelmann}, there is a primitive idempotent $f\in (\O' Gb)^{N_G(P,e)}$ such that ${\rm Br}_P(f)=e$; set $i=jf$, then $i$ is a source idempotent of $b$, and contained in $(\O'Gb)^P$. Multiplication by $f$ induces a unitary algebra homorphism
$$j\O'N_G(P)j\to i\O'Gi$$
which is split injective as a $(j\O'N_G(P)j,j\O'N_G(P)j)$-bimodule homomorphism (see \cite[Proposition 4.10]{FP} or \cite[Theorem 6.15.1]{Linckelmann}). So we can regard $i\O'Gi$ and $i\O'G$ as $j\O'N_G(P)j$-modules via the homomorphism.
}
\end{notation}

For finite groups $G,H$, an $(\O G,\O H)$-bimodule $M$ can be regarded as an $\O(G\times H)$-module (and vice versa) via $(g,h)m=gmh^{-1}$, where $g\in G$, $h\in H$ and $m\in M$. If $M$ is indecomposable as an $(\O G,\O H)$-bimodule, then $M$ is indecomposable as an $\O(G\times H)$-module, hence has a vertex (in $G\times H$) and a source. Denote by $(\O H)^\circ$ the opposite algebra of $\O H$ and by $\theta$ the $\O$-algebra isomorphism $\O H\cong(\O H)^\circ$ sending every $h\in H$ to $h^{-1}$. Let $b$, $c$ be blocks of $\O G$, $\O H$ respectively. Clearly $\theta(c)$ is a block of $\O H$. Then an $(\O Gb,\O Hc)$-bimodule $M$ can be regarded as an $\O(G\times H)$-module belongs to the block $b\otimes \theta(c)$ of $\O(G\times H)$, and vice versa. Here, we identify $b\otimes \theta(c)$ and its image under the $\O$-algebra isomorphism $\O G\otimes_\O \O H\cong \O(G\times H)$ sending $g\otimes h$ to $(g,h)$ for any $g\in G$ and $h\in H$. We will use the notation $\Delta G$ to denote the group $\{(g,g)|g\in G\}$.

\begin{lemma}\label{fOG}
Keep the notation of Notation \ref{notation}. Up to isomorphism, the $(\O'N_G(P,e)\hat{e}, \O'Gb)$-bimodule $f\O'G$ is the unique direct summand of $\hat{e}\O'Gb$ having $\Delta P$ as a vertex.
\end{lemma}

\noindent{\it Proof.} It is clear that the idempotent $b\hat{e}$ is in $(\O'Gb)^{N_G(P,e)}$. Since ${\rm Br}_P((b\hat{e})f)=e$, we have $(b\hat{e})f\neq 0$. Since $f$ is primitive in $(\O'Gb)^{N_G(P,e)}$ and $(b\hat{e})f=f(b\hat{e})$, we have $(b\hat{e})f=f$. Note that the endomorphism algebra of the $(\O'N_G(P,e),\O'Gb)$-bimodule $\O'Gb$ is isomorphic to $(\O'Gb)^{N_G(P,e)}$. So the $(\O'N_G(P,e)\hat{e}, \O'Gb)$-bimodule $f\O'G$ ($=f\O'Gb$) is an indecomposable direct summand of $\hat{e}\O'Gb$ ($=b\hat{e}\O'Gb$).

It follows from \cite[Proposition 6.7.4]{Linckelmann} that the bimodule $f\O'G$ has $\Delta P$ as a vertex. Note that $\hat{e}\O'Gb=f\O'G\oplus (b\hat{e}-f)\O'G$ and that ${\rm Br}_P(b\hat{e}-f)=e-e=0$.
By \cite[Theorem 3.2 (1)]{Broue}, any indecomposable direct summand of $(b\hat{e}-f)\O'G$ cannot have $\Delta P$ as a vertex.  $\hfill\square$

\medskip To prove the following lemma, we will use the concept of almost source algebras introduced by Linckelmann (\cite[\S 4]{Lin08}) and use the property that vertices and sources of indecomposable modules can be read off from almost source algebras (cf. \cite[Theorem 6.4.10]{Linckelmann}). Although the coefficient field in \cite{Lin08} is assumed to be algebraically closed, the definition of almost source algebras does not require the field to be algebraically closed (cf. \cite[Definition 6.4.3]{Linckelmann}). The ``algebraically closed" is also not used in the proof of \cite[Proposition 4.9 (i),(ii)]{Lin08}.

\begin{lemma}\label{summand}
Keep the notation of Notation \ref{notation}. Up to isomorphism, the $(\O'N_G(P)c,\O'Gb)$-bimodule $c\O'Gb$ has a unique indecomposable direct summand having $\Delta P$ as a vertex, which is isomorphic to $T':=\O'N_G(P)j\otimes_{j\O'N_G(P)j}i\O'G$.
\end{lemma}

\noindent{\it Proof.} By the standard Morita equivalences between $\O'N_G(P)c$ and $j\O'N_G(P)j$ as well as $\O'Gb$ and $i\O'Gi$, the functor (say $F_1$) sending an $(\O'N_G(P)c,\O'Gb)$-bimodule $M$ to the $(j\O'N_G(P)j,i\O'Gi)$-bimodule $jMi$ induces an equivalence between the category of $(\O'N_G(P)c,\O'Gb)$-bimodules and the category of $(j\O'N_G(P)j,i\O'Gi)$-bimodules.

By the standard Morita equivalences between $\O'N_G(P,e)\hat{e}$ and $j\O'N_G(P,e)j$ as well as $\O'Gb$ and $i\O'Gi$, the functor (say $F_2$) sending an $(\O'N_G(P,e)\hat{e},\O'Gb)$-bimodule $M$ to the $(j\O'N_G(P)j,i\O'Gi)$-bimodule $jMi$ induces an equivalence between the category of $(\O'N_G(P,e)\hat{e},\O'Gb)$-bimodules and the category of $(j\O'N_G(P,e)j,i\O'Gi)$-bimodules.

Noting that $j\O'N_G(P)j=j\O'N_G(P,e)j$, we can easily see that $F_1(T')\cong F_2(f\O'G)=i\O'Gi$ and that $F_1(c\O'Gb)=F_2(\hat{e}\O'Gb)=j\O'Gi$. Hence the functors $F_1$ and $F_2$ induce a bijection (say $\Phi$) between the isomorphism classes of indecomposable direct summands of the $(\O'N_G(P)c,\O'Gb)$-bimodule $c\O'Gb$ and those of the $(\O'N_G(P,e)\hat{e},\O'Gb)$-bimodule $\hat{e}\O'Gb$, which sends the isomorphism class of $T'$ to the isomorphism class of $f\O'G$.

It is easy to check that $N_{N_G(P)\times G}(\Delta P)=\Delta N_G(P)(1\times C_G(P))$ and $N_{N_G(P,e)\times G}(\Delta P)=\Delta N_G(P,e)(1\times C_G(P))$. So if a subgroup $R$ of $P\times P$ is $(N_G(P)\times G)$-conjugate to $\Delta P$, then $R$ is $(N_G(P,e)\times G)$-conjugate to $\Delta P$.

Let $\theta$ be the $\O'$-algebra isomorphism $\O'G\cong (\O'G)^\circ$ mapping any $g\in G$ to $g^{-1}$. Then $\theta(b)$ is also a block of $\O'G$, and $\theta(i)$ is a source idempotent of $\theta(b)$. By \cite[Proposition 4.9]{Lin08}, the interior $(P\times P)$-algebra $j\O'N_G(P)j\otimes_{\O'} \theta(i)\O'G\theta(i)$ is an almost source algebra of the block $c\otimes \theta(b)$ of $\O'(N_G(P)\times G)$, as well as the block $\hat{e}\otimes \theta(b)$ of $\O'(N_G(P,e)\times G)$. By \cite[Theorem 6.4.10]{Linckelmann} and the preceding paragraph, the bijective $\Phi$
preserves the property of ``having (resp. not having) $\Delta P$ as a vertex". Now the statement follows from Lemma \ref{fOG}. $\hfill\square$


\begin{lemma}\label{key}
Keep the notation of Notation \ref{notation}. Let $T':=\O'N_G(P)j\otimes_{j\O'N_G(P)j}i\O'G$. Assume that $\O'$ is finitely generated as an $\O$-module. Suppose that $b\in \O G$. Then there is an indecomposable direct summand $T$ of the $(\O N_G(P)c, \O Gb)$-bimodule $c\O Gb$, such that $T'\cong\O'\otimes_\O T$. 
\end{lemma}

\noindent{\it Proof.} This follows from Lemma \ref{summand} and \cite[Chapter III, Lemma 4.14]{Feit} (or \cite[Lemma 5.1]{Kessar_Linckelmann}).
 $\hfill\square$

\begin{lemma}\label{splendid}
Keep the notation of Notation \ref{notation}.  Assume that $k'$ is a splitting field for all subgroups of $G$, and that $\O'Gb$ is splendidly Morita equivalent to $\O'N_G(P)c$. Then the $(\O'N_G(P)c, \O'Gb)$-bimodule $T'$ in Lemma \ref{key} induces a splendid Morita equivalence between between $\O'N_G(P)c$ and $\O'Gb$.
\end{lemma}

\noindent{\it Proof.} Any two source algebras of a block are isomorphic as $\O$-algebras (see e.g. \cite[Theorem 6.4.4]{Linckelmann}). By \cite[Theorem 9.7.4]{Linckelmann}, we can see that the unitary algebra homorphism
$j\O'N_G(P)j\to i\O'Gi$ in Notation \ref{notation} is an isomorphism. Since the $(i\O'Gi,i\O'Gi)$-bimodule $i\O'Gi$ induces a Morita self-equivalence of $i\O'Gi$, using the isomorphism above, we can see that the $(j\O'N_G(P)j,i\O'Gi)$-bimodule $i\O'Gi$ induces a Morita equivalence between $j\O'N_G(P)j$ and $i\O'Gi$. By the standard Morita equivalences between $\O'N_G(P)c$ and $j\O'N_G(P)j$ as well as $\O'Gb$ and $i\O'Gi$, we obtain the statement.                     $\hfill\square$

\medskip\noindent{\it Proof of Thoerem \ref{3}.} This follows immediately from Lemma \ref{splendid}, Lemma \ref{key} and \cite[Proposition 4.5 (c)]{Kessar_Linckelmann}. $\hfill\square$

\medskip\noindent{\it Proof of Corollary \ref{corollary2}.} Assume that $\O$ is the minimal complete discrete valuation ring of $b$ in $\O'$ (see \cite[Definition 1.8]{Kessar_Linckelmann} for this conception). By \cite[Corollary 1.9]{Kessar_Linckelmann} and the remarks following it, we only need to show that $\O Gb$ and its Brauer corresponding block algebra are splendidly Morita equivalent. Denote by $\bar{b}$ the image of $b$ in $kG$. Since $p$-permutation modules of finite groups lift uniquely from $k$ to $\O$ up to isomorphism, it suffices to show that $kG\bar{b}$ is splendidly Morita equivalent to its Brauer corresponding block algebra. This follows from Theorem \ref{3}. $\hfill\square$

\begin{lemma}\label{Klein}
Keep the notation of Notation \ref{notation}. Denote by $b_0$ the principle block idempotent of $\O'A_5$.  Assume that $k'$ is a splitting field for all subgroups of $G$, and that $P$ is a Klein four group. If $i\O'Gi\cong \O'A_5b_0$, or $i\O'Gi\cong \O'A_4$, then $j\O'N_G(P)j\cong \O'A_4$; If $i\O'Gi\cong \O'P$, then $j\O'N_G(P)j\cong \O'P$.

\end{lemma}

\noindent{\it Proof.} Since we use \cite[Theorem 1.1]{CEKL} in this note, the statements can be deduced from \cite[\S 12.3]{Linckelmann}. $\hfill\square$

\medskip Let $A$ be an $\O$-algebra, and set $A':=\O'\otimes_{\O}A$, an $\O'$-algebra. Let $\Gamma:=\{\sigma\in {\rm Aut}(\O')~|~\sigma(u)=u, ~\forall~u\in \O\}$. For an $A'$-module
$U$ and a $\sigma\in \Gamma$, denote by ${}^\sigma U$ the $A'$-module which is equal to $U$ as a module over the subalgebra $1\otimes A$ of $A'$, such that $\lambda\otimes a$ acts on $U$ as $\sigma^{-1}(\lambda) \otimes a$ for all $a \in A$ and $\lambda\in \mathcal{O}'$.
An $A'$-module $M'$ is {\it $\Gamma$-stable} if ${}^\sigma M'\cong M'$ for all $\sigma\in \Gamma$. $M'$ is said to be {\it defined over $\O$}, if there is an $A$-module $M$ such that $M'\cong \O'\otimes_{\O} A$.

\medskip\noindent{\it Proof of Theorem \ref{3.1}.} We borrow the notation and the assumption in Theorem \ref{3.1} and Notation \ref{notation}. Since any block of a finite group algebra has a finite splitting field, we may assume that $k$ and $k'$ are finite. Since $p$-permutation modules of finite groups lift uniquely from $k$ to $\O$, up to isomorphism, it is sufficient to prove the statement for $\O=k$ and $\O'=k'$.

Let us firstly
consider the case that $G=A_5$, $b=b_0$, the principle block idempotent of $k'A_5$, and $P=\{(1),(12)(34),(13)(24),(14)(23)\}$. Then $N_G(P)=A_4$, $c=(1)$. Here we identify $A_4$ to the
subgroups of $A_5$ consisting by all elements which fix the letter 5. Identify $k'A_4$ with its image in $k'A_5b_0$ via multiplication by $b_0$.
Set $M':=k'Gb$, a $(k'N_G(P)c,k'Gb)$-bimodule. Let $\{T_1,T_2,T_3\}$ be a set of representatives of the isomorphism classes of simple $k'N_G(P)c$-modules, such that $T_1$ is the trivial module. For $i,j\in\{1,2,3\}$ such that $i\neq j$, let $T_j^i$ denote a uniserial $k'Gb$-module with composition series $T_i$, $T_j$, from top to bottom. For $1\leq i\leq 3$, denote by $Q_i$ a projective cover of $T_i$. By \cite[Corollary 12.2.10]{Linckelmann}, there is a set of representatives of isomorphism classes of simple $k'Gb$-modules $\{S_1,S_2,S_3\}$ with $S_1$ the trivial module, such that the restrictions of $S_1$, $S_2$, $S_3$ to $k'N_G(P)c$ are isomorphic to $T_1$, $T_3^2$, $T_2^3$, respectively. For $1\leq i\leq 3$, denote by $R_i$ a projective cover of the simple $k'Gb$-module $S_i^*:={\rm Hom}_{k'}(S_i,k')$. By \cite[Proposition 4.5.12]{Linckelmann}, a projective cover of $M'$ is isomorphic to $\oplus_{i=1}^3Q_i\otimes_{k'}R_i$. Set $Q'=\oplus_{i=2}^3Q_i\otimes_{k'}R_i$ and denote by $\pi':Q'\to M'$ the restriction to $Q'$ of a surjective bimodule homomorphism from this projective cover to $M'$. By \cite[\S 3]{Rickard} or \cite[Thereom 12.4.2]{Linckelmann}, the complex of $(k'N_G(P)c,k'Gb)$-bimodules of the form
$$X':=\cdots\to0\to Q'\xrightarrow{\pi'}M'\to 0\to \cdots$$
is a splendid Rickard complex. Here $M'$ lies in degree $0$ of $X'$.

By \cite[Proposition 4.5 (a)]{Kessar_Linckelmann}, to prove the theorem, it suffices to prove that there is a complex of $(kN_G(P)c,kGb)$-bimodules $X$ fulfilling $k'\otimes_kX\cong X'$. Set $M:= kGb$, a $(kN_G(P)c,kGb)$-bimodule. Then $M'\cong k'\otimes_k M$. Set $\Gamma:={\rm Gal}(k'/k)$. Note that $\Gamma$ preserves the properties of being simple, projective, trivial and indecomposable, hence permutes the sets $\{T_2, T_3\}$ and $\{S_2,S_3\}$, respectively, which in turn implies that $\Gamma$ permutes the sets $\{Q_2,Q_3\}$ and $\{R_2,R_3\}$, respectively. So $Q'$ is $\Gamma$-stable. It follows from \cite[Lemma 6.2 (c)]{Kessar_Linckelmann} that there is a projective $(kN_G(P)c,kGb)$-bimodule $Q$ such that $Q'\cong k'\otimes_k Q$. To show that $\pi'$ can be chosen to be of the form ${\rm Id}_{k'}\otimes \pi$ for some bimodule homomorphism $\pi : Q\to M$, consider a projective cover $\pi: Z\to M$. Then $k'\otimes_k Z\cong Z'$ yields the projective cover of $M'$ above. By \cite[Lemma 6.2 (c)]{Kessar_Linckelmann} and the Noether-Deuring Theorem (cf. \cite[page 139]{CR}), $Z$ has a direct summand isomorphic to $Q$. So we just need to restrict $\pi$ to $Q$, then $\pi$ is a desired map. This shows that $kN_G(P)c$ and $kGb$ are splendidly Rickard equivalent.

Now, let us consider the general case. If $ik'Gi$ is isomorphic to $k'P$ (resp. $k'A_4$), then $jk'N_G(P)j$ is isomorphism to $k'P$ (resp. $k'A_4$). It is easy to see that $k'Gb$ and $k'N_G(P)c$ are Morita equivalent via the bimodule
$$T':=k'N_G(P)j\otimes_{jk'N_G(P)j}ik'Gi\otimes_{ik'Gi}ik'G.$$
By Lemma \ref{key}, there is an indecomposable direct summand $T$ of the $(kN_G(P)c, kGb)$-bimodule $ckGb$, such that $T'\cong k'\otimes_k T$. By \cite[Proposition 4.5 (c)]{Kessar_Linckelmann}, $kGb$ and $kN_G(P)c$ are splendidly Morita equivalent.

Assume that $ik'Gi$ is isomorphic to $k'A_5b_0$, then $jk'N_G(P)j$ is isomorphic to $k'A_4$. We identify $ik'Gi$ with $k'A_5b_0$ and identify $jk'N_G(P)j$ with $k'A_4$ via the isomorphisms. Then, the $(k'A_5b_0,k'Gb)$-bimodule $k'A_5b_0\otimes_{ik'Gi}ik'G$ induces a splendid Morita equivalence between $k'A_5b_0$ and $k'Gb$; the $(k'N_G(P)c,k'A_4)$-bimodule $k'N_G(P)j\otimes_{jk'Hj} k'A_4$ induces a splendid Morita equivalence between $k'N_G(P)c$ and $k'A_4$. Set
$$\widetilde{Q}':=k'N_G(P)j\otimes_{jk'N_G(P)j}Q'\otimes_{ik'Gi}ik'G$$ and $$\widetilde{M}':=k'N_G(P)j\otimes_{jk'N_G(P)j}M'\otimes_{ik'Gi}ik'G.$$
Then the complex
$$Y':=\cdots\to0\to \widetilde{Q}'\xrightarrow{\widetilde{\pi}'} \widetilde{M}'\to 0\to \cdots$$
induces a splendid Rickard equivalence between $k'N_G(P)c$ and $k'Gb$. Here $\widetilde{M}'$ lies in degree $0$ of $Y'$. As in the first case, we are going to prove that there is a complex of $(kN_G(P)c,kGb)$-bimodules $Y$ satisfying $k'\otimes_kY\cong Y'$.

Note that $\widetilde{M}'=k'N_G(P)j\otimes_{jk'N_G(P)j}ik'Gi\otimes_{ik'Gi}ik'G$. By Lemma \ref{key}, $\widetilde{M}'$ is defined over $k$.
For $1\leq i\leq 3$, set $$\widetilde{Q}_i:=k'N_G(P)j\otimes_{jk'N_G(P)j}Q_i~~{\rm and}~~\widetilde{R}_i:=k'Gi\otimes_{ik'Gbi}R_i.$$
Since Morita equivalences preserve the properties of being simple, projective and compatible with projective covers, we can see the following things: $\widetilde{Q}_1,\widetilde{R}_1$ are projective covers of the trivial $k'N_G(P)c$-module and the trivial $k'Gb$-module, respectively; $\oplus_{i=1}^3\widetilde{Q}_i\otimes_{k'}\widetilde{R}_i$ is a projective cover of $\widetilde{M}'$; $\widetilde{Q}'=\oplus_{i=2}^3\widetilde{Q}_i\otimes_{k'}\widetilde{R}_i$,
and $\widetilde{\pi}':\widetilde{Q}'\to \widetilde{M}'$ is the restriction to $\widetilde{Q}'$ of the surjective bimodule homomorphism from the projective cover to $\widetilde{M}'$. Then the rest is analogous with the proof of the first case. $\hfill\square$

\medskip\noindent{\it Proof of Thoerem \ref{main}.} The proof is similar to the proof of \cite[Theorem 1.12]{Kessar_Linckelmann}. We only need to replace \cite[Theorem 1.10]{Kessar_Linckelmann} there with Theorem \ref{3.1}. $\hfill\square$


\bigskip\noindent\textbf{Acknowledgement.}\quad This note is written during the time of me as a visiting student at City, University of London in 2021. I wish to thank Professor Markus Linckelmann for suggesting me this topic, for plenty of patient answers to my questions and for many insightful suggestions.



\end{document}